\documentclass[leqno,draft]{article}

\begin{document}

\title{A few remarks about operator theory, topology, and analysis on
metric spaces}

\author{Stephen Semmes	\\
	Rice University}

\date{}

\maketitle

\begin{abstract}
Some basic facts about Fredholm indices are briefly reviewed, often
used in connection with Toeplitz and pseudodifferential operators, and
which may be relevant for operators associated to fractals.
\end{abstract}

	Let $\mathcal{H}$ be a complex infinite-dimensional separable
Hilbert space, which is of course unique up to isomorphism.  Let
$\mathcal{B}(\mathcal{H})$ be the Banach algebra of bounded linear
operators on $\mathcal{H}$, and let $\mathcal{C}(\mathcal{H})$ be the
closed ideal in $\mathcal{B}(\mathcal{H})$ consisting of compact
operators on $\mathcal{H}$.  The quotient $\mathcal{A}(\mathcal{H}) =
\mathcal{B}(\mathcal{H}) / \mathcal{C}(\mathcal{H})$ is a $C^*$
algebra known as the \emph{Calkin algebra}.

	A bounded linear operator $T$ on $\mathcal{H}$ is said to be
\emph{Fredholm} if the kernel of $T$ is finite-dimensional and $T$
maps $\mathcal{H}$ onto a closed linear subspace of $\mathcal{H}$ of
finite codimension.  Sometimes it is convenient to consider unbounded
linear operators too, such as differential operators.  It may also be
helpful to include operators between different Hilbert spaces.  One
way to deal with unbounded linear operators is to treat them as
bounded linear operators between different Hilbert spaces.

	If $T \in \mathcal{B}(\mathcal{H})$ is Fredholm, then there is
an $R \in \mathcal{B}(\mathcal{H})$ such that $R \, T - I$ and $T \, R
- I$ are finite rank operators, where $I$ is the identity operator on
$\mathcal{H}$.  Conversely, $T \in \mathcal{B}(\mathcal{H})$ is
Fredholm if there is an $R \in \mathcal{B}(\mathcal{H})$ such that $R
\, T - I$ and $T \, R - I$ are compact operators.  Equivalently, $T
\in \mathcal{B}(\mathcal{H})$ is Fredholm if and only if the
corresponding element $\widehat{T}$ of $\mathcal{A}(\mathcal{H})$ is
invertible.

	The \emph{index} of a Fredholm operator $T \in
\mathcal{B}(\mathcal{H})$ is defined to be the difference between the
dimension of the kernel of $T$ and the codimension of $T(\mathcal{H})$
in $\mathcal{H}$.  If $T, T' \in \mathcal{B}(\mathcal{H})$, $T$ is
Fredholm, and $T' - T$ is compact, then $T'$ is a Fredholm operator
with the same index as $T$.  If $T, T' \in \mathcal{B}(\mathcal{H})$,
$T$ is Fredholm, and $T'$ is sufficiently close to $T$ in the operator
norm topology, then $T'$ is again Fredholm and has the same index as
$T$.  The composition of two Fredholm operators is Fredholm, and the
index of the composition is equal to the sum of the indices.

	Many questions in analysis are concerned with invertability of
linear operators.  One might start by showing that an operator is
Fredholm, and then try to analyze the kernel and cokernel.  These may
be considered as nuisances to be minimized.  By contrast, from the
perspective of topology, nontrivial indices are an opportunity.

	As in finite dimensions, invertible operators on $\mathcal{H}$
form a connected open set in $\mathcal{B}(\mathcal{H})$.  Fredholm
operators form an open set that is not connected, with different
components distinguished by the index.

	The set of invertible operators on $\mathcal{H}$ is in fact
contractable \cite{kpr}.  The set of Fredholm operators is not, and
can instead be used as a classifying space for topological $K$-theory
\cite{a2}.

	Fredholm operators occur frequently in the study of integral
equations, where they may be given as compact perturbations of
multiples of the identity and hence have index $0$ automatically.
However, they may also have nontrivial kernels and cokernels, as in
the context of elliptic boundary value problems.

	A Toeplitz operator associated to a continuous complex-valued
function on the unit circle is Fredholm when the function is
nonvanishing.  The index of the Toeplitz operator is equal to the
winding number of the function around $0$ in the complex plane.

	There are very interesting ways in which to elaborate some
Fredholm operators and get additional indices.  One approach applies
to operators whose commutators with multiplication operators by
continuous functions are compact \cite{a1}.  This permits the operator
to be adjusted to act on spaces of sections of vector bundles.
Another approach deals with homomorphisms from algebras of continuous
functions on compact metric spaces into the Calkin algebra
\cite{b-d-f-1, b-d-f-2, b-d-f-3, d1}.  In this case there are also
operators corresponding to matrix-valued functions on the space.

	Techniques like these are very important for increasing the
information available from an operator or a collection of operators.
A basic construction may lead to a whole family of indices.

	It seems to me that there is a lot of room for development of
operator theory on fractals.  In particular, there are a lot of
interesting fractals with topological dimension $1$.  Of course, many
fundamental notions and results of analysis carry over to the broad
setting of spaces of homogeneous type \cite{c-w-1, c-w-2}.  Part of
the point would be to take advantage of extra structure when it is
present, even if it is not quite the same as the usual smooth case.
There can still be significant interactions with the underlying
geometry and topology.

        Extensions of $C^*$-algebras as in \cite{b-d-f-1, b-d-f-2,
b-d-f-2, d1} are especially appealing for spaces with topological
dimension $1$, where the higher-dimensional cohomology groups are
trivial.  It is also nice to be able to view a $1$-dimensional space
as some sort of boundary of a $2$-dimensional space.  This includes
fractals in the plane, and the spaces considered by Bourdon and Pajot
\cite{b, b-p-1, b-p-2, b-p-3, b-p-4}, for instance.  As in these
examples, the $2$-dimensional spaces may have singularities or other
complications.  At any rate, it would be interesting to have something
like complex-analytic metric space involved.

\end{document}